\documentclass{article}
\usepackage{latexsym,amsmath}

\long\def\@makecaption#1#2{\vskip 10pt \setbox\@tempboxa\hbox{#1. #2}  
   \ifdim \wd\@tempboxa >\hsize   
       #1. #2\par                 
     \else                        
       \hbox to\hsize{\hfil\box\@tempboxa\hfil}                    
   \fi}
\newtheorem{thm}{Theorem}[section]
\newtheorem{lem}[thm]{Lemma}

\newtheorem{cor}[thm]{Corollary}

\newcommand{\qed}{\hfill\mbox{$\Box$}}
\newenvironment{pf}{\medskip {\bf Proof.\ }}{\qed\medskip}
\newcommand{\boxx}{\Box}
\newcommand{\LHS}{left-hand side}
\newcommand{\RHS}{right-hand side}

\newcommand{\G}{\Gamma}
\newcommand{\Diam}{D}
\renewcommand{\l}{\ell_1}
\title{On equicut graphs}
\author{Michel Deza\thanks{CNRS and Ecole Normale Sup\'erieure, Paris, France}
\and Dmitrii V. Pasechnik\thanks{Dept. of Computer Science, Utrecht University,
The Netherlands}}
\begin{document}
\maketitle
\begin{abstract}
  The {\em size} $sz(\G)$ of an $\l$-graph $\G=(V,E)$ is the minimum of
  $n_f/t_f$ over all the possible $\l$-embeddings $f$ into
  $n_f$-dimensional hypercube with scale $t_f$.  The sum of distances
  between all the pairs of vertices of $\G$ is at most $sz(\G)\lceil
  v/2 \rceil\lfloor v/2 \rfloor$ ($v=|V|$).  The latter is an equality
  if and only if $\G$ is {\em equicut graph}, that is, $\G$ admits an 
  $\l$-embedding $f$ that for any $1\leq i\leq n_f$ satisfies $\sum_{x\in
    X} f(x)_i\in \{ \lceil v/2 \rceil, \lfloor v/2 \rfloor\}$ for any
  $x\in V$.  Basic properties of equicut graphs are investigated. A
  construction of equicut graphs from $\l$-graphs via a natural
  doubling construction is given.  It generalizes several well-known
  constructions of polytopes and distance-regular graphs.  Finally,
  large families of examples, mostly related to polytopes and
  distance-regular graphs, are presented.
\end{abstract}

\section{Background}
A finite simple graph $\G=(V,E)$ with $v=v(\G)$ vertices has the
natural graph distance denoted by $d(x,y)=d_\G(x,y)$.
The diameter of $\G$, that is, that maximum of $d(x,y)$ over
all the pairs of vertices $x$, $y$, is denoted by $\Diam(\G)$.
Then $\G=(V,E)$ is an
$\l$-{\em graph} if there exists a mapping $f: V\mapsto R^m$
with $d_\G(x,y)$ 
being equal to the $\l$-distance
between $f(x)$ and $f(y)$ for any two $x, y\in V$ or, equivalently, if
there exists a mapping $g: V\mapsto H^n$ (the latter being the
vertex-set of the $n$-cube) such that for any $x, y\in V$ we have
$td(x,y)=d_H(g(x),g(y))$ for some integer $t=t_f$. Recall that 
$d_H(g(x),g(y))$ is the Hamming distance, that is, the size of the
symmetric difference of the sets with characteristic vectors
$g(x)$ and $g(y)$. The number $t$ is
called the {\em scale} of the embedding $f$.
For more information on $\l$-graphs we refer the reader to 
\cite[Chapt.~21]{DeLa97}.

In fact, such $t$ should be 1 or even number that is at most $v-2$
(assuming $v\geq 4$).  Call the smallest such $t=t_f$ (if
any) the {\em scale} $t(\G)$ of $\G$.  Call the the minimum of $n/t$
the {\em size} $sz(\G)$ of $\l$-graph $\G$.  
(In general, $sz(\G)$ is
larger than minimal $m$, the dimension of the host $\l$-space into
which the shortest path metric of $\G$ embeds isometrically.)

$W(\G)$ denotes the sum of all $v(v-1)/2$ pairwise distances between
vertices of $\G$; chemists call it the {\em Wiener number}.
In \cite{Ple84} this number is called {\em transmission}.

A powerful tool to deal with $\l$-graphs is via linear programming.
Namely, consider the cone $Met(V)$ of {\em semimetrics} on $V$. 
The elements of $Met(V)$ are symmetric $v\times v$-matrices
with 0 diagonal that satisfy the triangle inequalities. A {\em metric} is
a semimetric that does not have 0's outside the main diagonal.
The distance matrix $Dist(\G)$ of a graph $\G$ is a point in $Met(V)$.

The set of {\em cuts} $S$ of $V$ consists of partitions of $V$ into two
parts.  Each $s\in S$ defines an extremal ray of the cone
$Met(V)$.
Namely, $s$ defines a
semimetric $\delta_s$ on $V$, so that $\delta_s(u,v)=0$ or $1$
depending on $u,v\in V$ being in same or different parts of $s$.  Then
one may ask whether 
$Dist(\G)$ lies in the cone $Cut(V)$ that is generated by 
the cuts, that is
\begin{equation}\label{l1emb0}
Dist(\G) = \sum_{s\in S}\lambda_s \delta_s
\end{equation}
for the appropriate choice of $\lambda_s$.
Hence $Dist(\G)\in Cut(V)$
implies that 
\begin{equation}\label{l1emb}
t\cdot Dist(\G) = \sum_{s\in S}\lambda_s \delta_s
\end{equation}
for an integer $t$ and an integer vector $\lambda=(\lambda_s\mid s\in S)$. 

The above gives another definition of the
$\l$-graphs: a graph $\G$ is $\l$-graph if and only if $Dist(\G)\in Cut(V)$.
Each equality (\ref{l1emb}) for a given $t$ and $\lambda$ defines an
$t$-scale embedding $f$
of $\G$ into the hypercube of dimension $\sum_{s\in S}\lambda_s$.

An advantage of using (\ref{l1emb}) is that is allows to classify
$\l$-embeddings of $\G$ up to equivalence: different solutions to (\ref{l1emb})
with nonnegative integer unknowns $\lambda$ such that $g.c.d.(t,\lambda_i)=1$
correspond to different embeddings. If such a solution is unique
$\G$ is called {\em rigid}. 
Any rigid $\l$-graph has scale 1 or 2; any $\l$-graph
of scale 1 is rigid. If an $\l$-graph does not contain $K_4$ then it is rigid.
Cf. \cite[Sect.~21.4]{DeLa97} for further information.

It is easy to check the following reformulation of \cite[(4.3.6)]{DeLa97}.
\begin{lem}\label{i3}
$W(\G)/(\lceil v/2 \rceil\lfloor v/2 \rfloor) \leq
sz(\G)\leq W(\G)/(v-1).$
\end{lem}
\begin{pf} 
Let the ($t$-scale) 
embedding $f$ be given by the matrix $F=(F_{iu})=(f(u)_i)$.
The columns of $F$ can be regarded as characteristic vectors of
one of the two parts of cuts $s\in S$.
Using (\ref{l1emb}) and 
summing up the $d_\G (u,v)$'s, we obtain
\begin{eqnarray*}
W & = &W(\G)=\sum_{u,v\in V} d_\G (u,v) = 
\sum_{u,v\in V} \sum_{s\in S}\lambda_s \delta_s(u,v) \\ 
& = & \sum_{s\in S} \sum_{u,v\in V} \lambda_s \delta_s(u,v) = 
\sum_{s\in S} \lambda_s |s|(v-|s|) \label{eq:00-0} \\
& \leq &
\lceil v/2 \rceil\lfloor v/2\rfloor\sum_{s\in S}\lambda_s,
\end{eqnarray*}
by taking the maximum over $|s|$ in (\ref{eq:00-0}).
On the other hand, by taking the minimum there, we obtain
$W\geq (v-1)\sum_{s\in S}\lambda_s.$
\end{pf}

When the equality holds in the \LHS\ of the inequality in
Lemma~\ref{i3}, we say that $\G$ is an {\em equicut graph}.
This means that for such a graph in the formulae (\ref{l1emb})
every $s\in S$ satisfies $\delta_s\not=0$ if and only if $s$ partitions $V$ into
parts of size $\lceil v/2 \rceil$ and $\lfloor v/2 \rfloor$.

On the other hand, the equality on the other side only happens in
very special case.
\begin{lem}\label{star}
$sz(\G) = W(\G)/(v-1)$ if and only if $\G$ is the star $K_{1,v-1}$.
\end{lem}
\begin{pf}
Without loss in generality, the first row $f(u)$ of the matrix $F$ defined by
$f$ consists of 0's. Moreover, all the cuts $s$ that appear in $F$ are
partitions of $V$ into part of size one and of size $v-1$.
It follows that the any geodesic between $x,y \in V$ must contain
$u$. Hence $\G$ is a star.
\end{pf}

Further, the size $sz(\G)$ has the following properties formulated
in terms of $v$.
\begin{thm} \label{prop1}
Let $\G$ be an $\l$-graph.
\begin{itemize}
\item[$(i)$] $sz(\G)=2-1/(\lceil v/2 \rceil)$ if and only if $\G=K_v$;
\item[$(ii)$] $sz(\G)=2$ if and only if $\G$ is a non-complete subgraph of a 
  Cocktail-Party graph;
\item[$(iii)$] $sz(\G) = v-1$ if and only if $\G$ is a tree;
\item[$(iv)$] $2<sz(\G)<v-1$ otherwise.
\end{itemize}
\end{thm}
\begin{pf}
By taking all the equicuts on $V$, one obtains a realization of
$K_v$ with $n/t=2-1/(\lceil v/2 \rceil)$. 
By using the \LHS\ of the inequality
in Lemma~\ref{i3} and the fact that all distances in $K_v$ are 1,
we obtain that $n/t=sz(\G)$, as required in $(i)$.
\medskip

The size of the Cocktail Party graph is 2, as it is easy to construct
such a realization (see for instance, Theorem~\ref{lem3} below).  Then
observe that $sz(\G)<2$ implies that $\G$ is complete, as a two-path
from $x$ to $y$ the vertex $x$ is realized, without loss of
generality, by the zero vector, and the vertex $y$ is realized by a
vector with $2t$ ones.  Hence any non-complete subgraph subgraph of
the Cocktail Party graph has size 2, and $(ii)$ is proved.  
\medskip

To prove the last two claims, consider the decomposition (\ref{l1emb0}).
Choose an $s^*\in S$ with $\lambda_{s^*}>0$; 
it defines the partition $V^+(s^*)\cup V^-(s^*) = V$.
It follows that the subgraph
of $\G$ induced on $V^+=V^+(s^*)$ is isometric. 
Indeed, consider a shortest path between 
$x^+$ and $y^+$. Let $k\in V$ lie on this path. Then
\begin{eqnarray*}
d_\G(x^+,k) & + & d_\G(k,y^+)-d_\G(x^+,y^+)  =  \\
& = & \sum_{s\in S}\lambda_s(\delta_s(x^+,k)+\delta_s(k,y^+)
 -\delta_s(x^+,y^+)) =  0.
\end{eqnarray*}
As $\lambda\geq 0$, we have 
$\delta_{s^*}(x^+,k)+\delta_{s^*}(k,y^+)=\delta_{s^*}(x^+,y^+)$.

By definition of cut, \RHS\ of the above is 0. Hence both the summands in
\LHS\ are 0, and $k\in V^+$. as required. Similarly, the subgraph
induced on $V^-(s^*)$ is isometric.

For an edge $(x,y)\in V^+\times V^-$ we have $d_\G(x,y)=1\geq\lambda_s$.
(Such an edge exists, as $\G$ is connected.)
Hence $\lambda\leq 1$.

Let ${\cal S}\subseteq S$ denote the set of $s\in S$ satisfying $\lambda_s>0$.
Now we proceed by induction on $v$.
Let $s'\in {\cal S}$. Hence the subgraph $\Delta^+$ of $\G$ induced on $V^+(s')$ is
isometric. Thus
\[
Dist(\Delta^+)=\sum_{s\in{\cal S}-\{s'\},V^+(s)\cap V^+(s')\not=\emptyset}
\lambda_s\delta_{V^+(s)\cap V^+(s')}. 
\]
By the induction assumption, 
\[
X^+ := \sum_{s\in{\cal S}-\{s'\},V^+(s)\cap V^+(s')\not=\emptyset}
\lambda_s\leq |V^+(s')|-1.
\]
Similarly, we obtain
\[
X^- := \sum_{s\in{\cal S}-\{s'\},V^+(s)\cap V^-(s')\not=\emptyset}
\lambda_s\leq |V^-(s')|-1.
\]
Hence 
\begin{eqnarray}\label{ineq1}
sz(\G)=\sum_{s\in\S}\lambda_s&=&
\lambda_{s'} + X^+ +X^- \leq 1+|V^+(s')|-1+|V^-(s')|-1=\\
&=&v-1.\nonumber
\end{eqnarray}
Finally, if  $sz(\G)=v-1$  then (\ref{ineq1}) turns into
an equality, and $\lambda_{s'}=1$. The latter and the induction implies
that $\G$ is a tree.
\end{pf}

For a tree $\G$, one can moreover show the following.
\begin{lem}
$W(K_{1,v-1})=(v-1)^2$, $W(P_v)=v(v-1)(v+1)/6$.
For any tree $\G$ that is neither the star $K_{1,v-1}$ nor the path
$P_v$, the value $W(\G)$ lies strictly within the extremities above.
\end{lem}

\section{Generalities on equicut graphs}
Recall that a connected 
graph $\G=(V,E)$ is called {\em $2$-connected} (or
{\em $2$-vertex-connected}) if it remains connected after deletion of any
vertex. 

\begin{lem}\label{lem20}
An equicut graph $\G$ with a least $4$ vertices is $2$-connected.
\end{lem}
\begin{pf}
Assume that $\G$ is not $2$-connected. Then there exists a vertex $x\in V$
so that the subgraph induced on $V-\{x\}$ has two connected components,
$V_1$ and $V_2$.  We can assume, without loss of generality, that
$x$ is represented by all-0 vector $F(x)$ in the equicut realization $F$ of $\G$.
Then for any $v_i\in V_i$ ($i=1,2$) one has
\[
d_\G(v_1,v_2)=d_\G(v_1,x)+d_\G(x,v_2)=|v_1|+|v_2|,
\]
where $|v_i|$ denotes the number of 1's in the vector $F(v_i)$.
Hence the vectors $F(v_1)$ and $F(v_2)$ have disjoint supports, for any
$v_1\in V_1$ and $v_2\in V_2$. That is, any column of $F$ either has
1's in $V_1$, or in $V_2$, but never in both of them.
As $F$ is an equicut, we obtain $|V_i|\leq 1$, as required.
\end{pf}

This implies
\begin{cor}\label{lem2}
For any equicut graph $\G$ with $v\geq 4$ vertices, we have
$$2-1/(\lceil v/2 \rceil)\leq sz(\G)\leq v/2 $$
with equality on the left if and only if $\G=K_v$ and on the right if and only if $\G=C_v$.
\end{cor}
\begin{pf}
The left hand side of the inequality follows from Lemma~\ref{prop1}.
Then, to see that $K_v$ is equicut, it suffices to observe that the embedding
$f$ given by the matrix $F$ with the columns being the cuts with the smallest
part of size $\lfloor v/2 \rfloor$ is an equicut embedding of $K_v$.

To see that $C_v$ is equicut, observe that $f$ given by the matrix $F$
obtained by the cyclic shifts of the row of the form
($\lceil v/2 \rceil$~times~1,$\lfloor v/2 \rfloor$~times~0)
produces an equicut realization. 

To show the last claim of the lemma and the right hand inequality, 
is suffices to use Lemma~\ref{lem20} and  the result from
\cite{Ple84} that the Wiener number $W(\G)$ of a $2$-connected graph on $v$ vertices 
is maximal for the circuit $C_v$ (see \cite[Th.~5]{Ple84}).
\end{pf}

The condition $v\geq 4$ is necessary in the statements above.  
Indeed, $sz(P_3)=2>sz(C_3)=3/2$. Note
that $P_2$ and $P_3$ are the only equicut trees. Also,
$W(C_5)=15<W(P_{123452})$ and $sz(C_5)=5/2<sz(P_{123452})$, where
$P_{123452}$ denotes the circle on $2\dots 5$ with an extra edge attached
to the vertex $2$.
   
\paragraph{Remark.} $\G$ is equicut graph if there is a realization with
the binary matrix $F$ with the column sums $\lceil v/2 \rceil$ or 
$\lfloor v/2 \rfloor$.  If, instead
of this condition, we asked that any {\em row} of $F$ has exactly $k$
1's, then we obtain other special $\l$-graph. Namely, one which embeds
isometrically up to scale $t$, into the Johnson graph $J(v,k)$.
It is observed in \cite{Shppers} that such graphs can be 
recognized in polynomial time
using the algorithm in \cite{Shp93}; the ``atom graph'' of
$\G$ is bipartite.

However, we are not aware of any similar characterization of the equicut 
graphs.
  
\section{Doublings of $\l$-graphs}
Call an equicut $\l$-graph an {\em antipodal doubling} if its realization 
(i.e. above $(0,1)$-matrix)
has form $F=\left(\begin{smallmatrix} A & O\\ J-A & J'\end{smallmatrix}\right)$, 
where $A$ is a $v/2\times n'$ 
(0,1)-matrix,  $J,J'$ are matrices  
consisting of 1's only and $O$ is $v/2 \times (n-n')$ matrix consisting of 0's.

If moreover the matrix $A$ is a realization, with the same scale $t$,
of a graph $\G'$, then it is straightforward to check that $J'$ has
$t(\Diam(\G')+1)-n$ columns.  Note that Double Odd graph $DO_{2s+1}$
(see \cite{BCN89}) with $s\geq 3$ is an example of antipodal doubling
with the matrix $A$ not corresponding to the realization of a graph
$\G'$, for any decomposition of $F$ of the above form.

\paragraph{Remark.} An antipodal
doubling is exactly an $\l$-graph that admits an antipodal isomorphism, i.e.  it
has a central symmetry (for any vertex, there is exactly one other on
the distance equal to the diameter) and the mapping of all vertices
into their antipodes is an isomorphism. Antipodal extensions of arbitrary 
$\l$-metrics was considered in \cite[Sect.~7.2]{DeLa97}.
\medskip

In order to investigate when one can construct an $\l$-graph from
an $\l$-graph via the antipodal doubling (see Theorem~\ref{lem3} below), let us introduce the
following definition.  For a graph $\G=(V,E)$, define its {\em
  diametral doubling} as the graph $\boxx\G$ with the vertex set
$V^+\cup V^-$ (where $V^*$ is a copy of $V$) and the adjacency as
follows: $x^\mu$ is adjacent to $y^\epsilon$ if $\mu=\epsilon$ and
$(x,y)\in E$, or if $\mu\not=\epsilon$ and $d_\G(x,y)=\Diam(\G)$.

\begin{lem} \label{cond-i}
The subgraphs of $\boxx\G$ induced on $V^*$ are {\em isometric}
to $\G$ if and only if 
\begin{eqnarray}\label{eq:cond-i}
d_\G(x,y) & \leq & 2 + d_\G(z_1,z_2) \mbox{ for any }
x,y,z_1, z_2\in V\\
& & \mbox{ satisfying } d_\G(x,z_1)=d_\G(y,z_2)=D(\G).\nonumber
\end{eqnarray}
\end{lem}
\begin{pf}
Clearly, the condition of the lemma is necessary for the subgraphs
induced on $V^*$ to be geodetic. Indeed, otherwise there can be a shorter
path between, say, $x^+$ and $y^+$ that passes through a vertex in $V^-$.

To prove sufficiency, observe that any geodetic between $p^+$ and $q^+$
either lies wholly in $V^+$, or contains a path of the form
$x^+,z_1^-,\dots,z_2^-,y^+$, with the (sub)path from $z_1^-$ to $z_2^-$ lying
wholly in $V^-$. In the former case there is nothing to prove.
In the latter case, by definition of $\boxx\G$, one has that
$d_\G(x,z_1)=d_\G(y,z_2)=\Diam(\G)$. By the assumption,
there is a path of length at most $2+d_\G(z_1,z_2)$ between
$x$ and $y$. Thus the geodetic can be replaced by another geodetic
with the path from $x^+$ to $y^+$ lying wholly in $V^+$.
\end{pf}

\begin{lem}\label{cond-ii}
Let $\G$ satisfy the condition of Lemma~\ref{cond-i}. Then
\begin{equation}\label{eq:cond-ii}
d_{\boxx\G}(x^+,y^-)=\Diam(\G)+1-d_\G(x,y)
\end{equation} 
if and only if
any geodetic in $\G$ lies on a geodetic of length $\Diam(\G)$.
\end{lem}
\begin{pf}
We begin by proving that a geodetic (in $\G$) between 
$x$ and $y$ can be extended to a geodetic of length $\Diam(\G)$.
By Lemma~\ref{cond-i}, there is a geodetic from $x^+$ to $y^-$ of the form
$x^+,\dots,z_1^+,z_2^-,\dots,y^-$ such that the path from $x^+$ to $z_1^+$
(respectively, from $y^-$ to $z_2^-$) 
lies wholly in $V^+$ (respectively, in $V^-$). 

We claim that there is a geodetic between $z_1$ and $z_2$ (by
definition of $\boxx\G$ they are at distance $D(\G)$) that
contains $x$ and $y$.
Clearly, 
\[
d_{\boxx\G}(x^+,y^-)=d_\G(x,z_1)+1+d_\G(z_2,y).
\]
Using (\ref{eq:cond-ii}), we derive
\[
D(\G)+1-d_\G(x,y)=d_\G(x,z_1)+1+d_\G(z_2,y).
\]
Hence
\[
D(\G) = (d_\G(z_1,z_2) = )\ d_\G(x,z_1)+d_\G(z_2,y)+d_\G(x,y),
\]
as claimed.

Now, suppose that a geodetic between $x$ and $y$ can be extended to a
geodetic of length $D(\G)$, say, between $z_1$ and $z_2$.
By reversing the argument above, we derive (\ref{eq:cond-ii}).
\end{pf}

Certain properties of $\G$ are inherited by $\boxx\G$.
\begin{lem}\label{proper-doublings}
Let $\G=(V,E)$ be a graph satisfying 
$(\ref{eq:cond-i})$ and $(\ref{eq:cond-ii})$.
Then $\boxx\G$ satisfies $D(\boxx\G)=D(\G)+1$, 
$(\ref{eq:cond-i})$ and $(\ref{eq:cond-ii})$.
\end{lem}
\begin{pf}
As $(\ref{eq:cond-i})$ and $(\ref{eq:cond-ii})$ hold for $\G$, we have
that each $p^+\in V^+$ has exactly one vertex $p$ at the maximal
distance $D(\boxx\G)=D(\G)+1$. Now checking $(\ref{eq:cond-i})$ is
straightforward. 

To check $(\ref{eq:cond-ii})$, take a geodetic
between $x^+$ and $y^-$. Note that
\[
d_{\boxx\G}(x^+,y^-)=D(\G)+1 - d_\G(x,y).
\]
Hence 
\begin{equation}\label{eq:011}
d_{\boxx\G}(x^+,y^-)=d_{\boxx\G}(y^+,y^-) - d_{\boxx\G}(x^+,y^+),
\end{equation}
and the geodetic in question lies on the path from $y^+$ to $y^-$, as
required.
For a geodetic between $x^+$ and $y^+$ observe that (\ref{eq:011}) also
proves that it is extendable to the path from $y^+$ to $y^-$.
\end{pf}

If $D(\G)=2$ then $\G$ satisfies $(\ref{eq:cond-i})$.
Moreover, $\G\not=K_v$ 
satisfies $(\ref{eq:cond-ii})$, unless it has an edge $(x,y)$ 
with $\G(x)=\G(y)$. In particular, any strongly regular graph
satisfies $(\ref{eq:cond-i})$ and $(\ref{eq:cond-ii})$.
Note that $\l$-graphs form a rather small sub-family of strongly
regular graphs.

\paragraph{Switchings.}
Not always $\Delta=\boxx\G$ is uniquely defining the graph $\G$ it was
constructed from; it can be that $\Delta=\boxx\G'$ for $\G'\not\cong\G$.
See Section~\ref{sect:examples} for many examples of this situation.

$\G$ and $\G'$ are related by the following graph operation.
The {\em diametral switching} of a graph $\G$ with
respect to $S\subset V(\G)$ is a graph $\G'$ that is obtained from $\G$ by
retaining the edges that lie within $S\times S\cup (V-S)\times (V-S)$
and replacing the set of edges from $S\times (V-S)$ with the set
$\{(x,y)\in S\times (V-S)\mid d_\G(x,y)=D(\G)\}$.  

Note that {\em Seidel switching} (see \cite{BCN89}) is an operation that
coincides with the diametral switching for graphs of diameter 2.

\begin{thm}\label{lem3}
Let $\G=(V,E)$ be a $\l$-graph. 
Then $\boxx\G$ is an $\l$-graph if $\G$ satisfies
$(\ref{eq:cond-i})$, $(\ref{eq:cond-ii})$ and
\begin{equation}\label{eq:size}
sz(\G)\leq\Diam(\G)+1,
\end{equation}
Moreover, if $\boxx\G$ is an $\l$-graph then
\begin{itemize}
\item $D(\boxx\G)=D(\G)+1$, $\boxx\G$ satisfies
$(\ref{eq:cond-i})$, $(\ref{eq:cond-ii})$ and $(\ref{eq:size})$ with {\em equality}, and
\item All the $\l$-realizations of $\boxx\G$ are equicut, 
of the form, up to permutation of rows and columns, 
and taking complements of columns, 
$\left(\begin{smallmatrix} A &O\\ J-A &J' \end{smallmatrix}\right)$ 
with $t(D(\G)+1)$ columns, so that
$A$ is an $\l$-realization of $\G$ with scale $t$.
\end{itemize}
\end{thm}
\begin{pf}
It is obvious from $sz(\G)\leq\Diam(\G)+1$ that 
there exists a ($v\times n$, say) matrix $A$ that determines an 
$\l$-embedding of $\G$
such that $n/t_A\leq \Diam(\G)+1$. Let $\lambda = t_A(\Diam(\G)+1) - n$ and
$B=\left(\begin{smallmatrix} A &O\cr J-A &J' \end{smallmatrix}\right)$, 
where $J'$ has $\lambda$ columns.
Identify the first $v$ rows of $B$ with $V^+$ and the last $v$ rows of
$B$ with $V^-$.

By construction, the Hamming distances between rows of $B$ are equal,
up to scale $t=t_A$, to the corresponding distances in $\boxx\G$.
Hence $\boxx\G$ is an $\l$-graph. In view of Lemma~\ref{proper-doublings},
$(\ref{eq:cond-i})$ and $(\ref{eq:cond-ii})$ hold for $\boxx\G$.

$\boxx\G$ satisfies $(\ref{eq:size})$ with equality, 
as by construction, $sz(\boxx\G)=D(\G)+1$.

Taking an $\l$-realization $B'$ of $\boxx\G$, we can reorder the rows of $B'$ 
so that, as in $B$, the first half of the rows corresponds to $V^+$, and the
second half corresponds to $V^-$. By reordering the columns and possibly taking
their complements, we reduce $B'$ to the desired form. The rest of the claim
follows from the fact that $\G$ induced on $V^+$ is isometric in
$\boxx\G$.
\end{pf}

\paragraph{Remarks.} 
$K_4-P_3$, $K_4-P_2$, the Dynkin diagram $E_6$
are examples of $\l$-graphs satisfying (\ref{eq:cond-i}), (\ref{eq:size}),
but not (\ref{eq:cond-ii}).
{\em Affine} Dynkin diagram $\tilde E_6$
An example of a graph that does not satisfy (\ref{eq:cond-i}) and (\ref{eq:size}),
but does satisfy (\ref{eq:cond-ii}).

Note that any $\G$ with $D(\G)=2$ satisfies (\ref{eq:cond-i}) and
(\ref{eq:cond-ii}). Certainly, not all of them are $\l$-graphs, 
for instance, $K_{2,3}$. Also, not all the $\l$-graphs of diameter
2 satisfy (\ref{eq:size}), for instance, $K_{1,4}$.

In general, for any $\l$-graph $\G=(V,E)$ with
$|V|\geq 4$ one has
\begin{equation}\label{eq:DG}
\Diam(\G)\leq sz(\G)\leq \Diam(\G)+|V|-3.
\end{equation}
The equality at the \RHS\ of (\ref{eq:DG}) holds if and only if $\G$ is a star, as
can be seen by applying Theorem~\ref{prop1}

Theorem~\ref{lem3} generalizes to arbitrary $\l$-graphs $\G$ the situation for
the Cocktail Party graph $K_{n\times 2}$ considered in \cite[Sect.~7.4]{DeLa97},
see also second part of Lemma~\ref{prop1}. It implies that the minimal
scale of an $\l$-embedding of $\boxx\G$ equals
the minimal $t>0$ so that the metric $t\cdot Dist(\G)$ embeds isometrically in
$H^{t(D(\G)+1)}$. In particular case of $\G=K_{4a}$ 
Lemma~7.4.6 of \cite{DeLa97}
states that such a minimal $t\geq 2a$, with the equality holding if and only if
there exists a Hadamard matrix of order $4a$.

\section{Examples}\label{sect:examples}
In this section we list many examples in no particular order. An equicut $\l$-graph $\G$ is called
a {\em doubling} if it is an antipodal doubling; moreover, if it is obtained from a graph $\Delta$ as
in Theorem~\ref{lem3}, we give such a representation $\G=\boxx \Delta$.

\paragraph{All equicut graphs with at most 6 vertices} are $C_v$ 
$(2\leq v\leq 6)$, $K_4$,
$K_5$, $K_6$, $P_3$, 4-wheel and the octahedron.  Among those 11
graphs only $K_4$, $K_5$, $K_6$ are not rigid 
and only $C_4$, $C_6$ and the octahedron are
doublings. 

\paragraph{Cartesian product} 
$\G\times\G'$ of two $\l$-graphs $\G$ and $\G'$ is
$\l$-graph if and only if both $\G$, $\G'$ are (see \cite[Sect.~7.5]{DeLa97});
the scale will be least common multiple of their scales and the size
will be the sum of their sizes. Moreover, $\G\times\G'$ is rigid if and only if they are.

\begin{lem}\label{prod-eq}
If $\G$, $\G'$ have even number of vertices each, then $\G\times\G'$
is an equicut graph if and only if they are. 
\end{lem}
\begin{pf}
\newcommand{\fblk}[1]
{\begin{smallmatrix}F_{#1}\\ \dots\\ F_{#1}\end{smallmatrix}}
\newcommand{\FFprod}{\left(\begin{smallmatrix}
    \fblk{1} & F' \\
    \fblk{2} & F'\\
    \dots    & \dots \\
    \dots    & \dots \\
    \fblk{v} & F'
    \end{smallmatrix}\right)}
Let $F$ and $F'$ be $\l$-realizations of $\G$ and $\G'$, respectively.
Then $B=\FFprod$, 
where $F_i$'s are rows of $F$, 
is an $\l$-realization of $\G\times\G'$.
Now let $F$ and $F'$ be equicut realizations.
It is straightforward that $B$ is an equicut realization when $\G$ and $\G'$ both
have even number of vertices.
\end{pf}

A {\em Doob} graph (the Cartesian product of a number of copies
of {\em Shrikhande graph} and a number of copies of $K_4$, see e.g.
\cite[p.27]{BCN89}) is an example of an equicut graph obtained via
Lemma~\ref{prod-eq}. It is a (non-rigid and non-doubling) $\l$-graph
of scale 2.

\paragraph{Distance-regular graphs.}
Here we freely use notation from \cite{BCN89}. Also, a significant
use is made of \cite{KoSh94}.
The Petersen graph and the Shrikhande graph are both
equicut graphs of scale 2 and size 3; both are rigid and are not
doublings.  The Double Odd graph $DO_{2s+1}$ is an equicut graph of scale 1 and
size $2s+1$.  The half-$s$-cube $\frac{1}{2} H(s,2)$ 
is an equicut
graph of scale 2 and size $s$.  It is not rigid only for $s=3,4$; it
is a doubling if and only if $s$ is even.  The Johnson graph $J(2s,s)$ is an
non-rigid doubling of scale 2 and size $s$.  

Further, the following are distance-regular  equicut graphs.
\begin{enumerate}
\item All Taylor $\l$-graphs - $\frac{1}{2} H(6,2)$, $J(6,3)$, $C_6$, 
$H(3,2)$, the Icosahedron. They are all doublings of diameter 3 and size 3
that can be constructed using Theorem~\ref{lem3} above.
\item All strongly regular $\l$-graphs,
except $J(s,2)$ ($s\geq 5$) and
  grids $H(2,s)$ with $s$ odd. 
That is, $C_5$, Petersen, $\frac{1}{2} H(5,2)$, Shrikhande, 
$H(2,s)$ with $s$ even, and $K_{s\times 2}$.
\item Among distance regular graphs 
$\G$ with $\Diam(\G)>2$ and $\mu>1$: $\frac{1}{2}H(s,2)$ with $s>5$, 
$H(s,d)$, $J(s,t)$ with $t>2$, the Icosahedron, and Doob graphs.
\item All Coxeter $\l$-graphs except $J(s,t)$ with $t<s/2$: 
$J(2s,s)$, the Icosahedron, the Dodecahedron, $K_{s\times 2}$, 
$\frac{1}{2} H(s,2)$, $H(s,2)$, $C_s$ ($s\geq5$).
\item All cubic distance-regular $\l$-graphs: $K_4$, Petersen, $H(3,2)$, 
$DO_5$, the Dodecahedron.
\end{enumerate}
Also all amply regular $\l$-graphs with $\mu>1$ are equicut graphs.
Yet another example is given by 
the 12-vertex co-edge regular subgraph of the Clebsh graph
$\frac{1}{2} H(5,2)$, see \cite[Sect.~3.11, p. 104]{BCN89}; it is an equicut graph
of size $5/2$, scale 2, non-doubling.

\paragraph{Some equicut graphs which are doublings of $\l$-graphs}
(see some in Section 7.2 of \cite{DeLa97}):
\begin{enumerate}
\item 
$C_{2s}=\boxx P_s$.
\item  $K_{s\times 2}=\boxx K_s$.
\item $H(s,2)=\boxx H(s-1,2)$.
\item $J(2s,s)=\boxx J(2s-1,s)$. \label{ex:J}
\item $\frac{1}{2} H(2s,2)= \boxx \frac{1}{2} H(2s-1,2)$.
\label{ex:halfH}
\item $Prism_{2s}=\boxx C_{2s}=\boxx L_{2s}$; here $L_{2s}$ is the
``ladder''on $2s$ vertices. $L_{2s}$ is the diametral switching of $C_{2s}$
with respect to the vertex set of a path $P_s\subset C_{2s}$.
\item $APrism_{2s+1}=\boxx C_{2s+1}$.
\item The Dodecahedron is the doubling of the 9-circle $(1,\dots,9)$ with 
extra vertex connected to the vertices 3, 6, 9 of the circle.
\item The Icosahedron is the doubling of 5-wheel;
as well, it is the doubling of the graph obtained from
hexagon $(1,\dots,6)$ by adding edges $(2i,2j)$ for $0<i<j<4$.  \label{tay9}
\item $J(6,3)$ is the doubling of Petersen graph, in addition to 
\ref{ex:J}.\label{tay10}
\item $ \frac{1}{2}H(6,2)$ is the doubling of Shrikhande graph, 
and the doubling
of $H(2,4)$ (more precisely, of its realization in half-$6$-cube) 
in addition to \ref{ex:halfH}. \label{tay11}
\end{enumerate}
In \ref{tay9}, \ref{tay10}, \ref{tay11} we have (diametral)
switching-equivalent graphs $\G$ such that $\boxx\G$ is a Taylor
$\l$-graph; see remark preceding Theorem~\ref{lem3}.  
This situation in general is well-known; for instance,
the Gosset graph, a Taylor graph that is not an $\l$-graph, 
can be obtained as 
the diametral doubling of one of 5 non-isomorphic, but switching-equivalent,
graphs. For definitions
and discussion of this situation in more general setting, see
e.g. \cite[pp.103-105]{BCN89}.

\paragraph{Equicut polytopes.}
The skeletons of many nice polytopes are equicut graphs.  Below we
list several such examples. We follow the terminology from 
\cite{Joh66,Cox73}.

All five Platonic solids have equicut skeletons; all, except the Tetrahedron, 
are rigid. All but the cube (of scale 1), have scale 2. The sizes for the
Tetrahedron, the Octahedron, the Cube, the Icosahedron and the Dodecahedron
are $3/2$, 2, 3, 3 and 5. 
 
The skeleton of any {\em zonotope} (see e.g. \cite{Zieg} for a definition of 
zonotope) is a doubling and it has scale 1, so it is rigid. 

Among the Archimedean $\l$-polytopes
\begin{enumerate}
\item all zonohedra (i.e. 3-dimensional zonotopes) are as follows:
  the truncated Octahedron, the truncated Cuboctahedron, the truncated
  Icosidodecahedron and $Prism_{2s}$ ($s>2$) with sizes 6, 9, 15 and $s+1$, 
  respectively; 
\item all the other doublings are as follows: the Rhombicuboctahedron, the 
Rhombicosidodecahedron and $Aprism_{2s+1}$ $(s>1)$ with
scale 2 and sizes 5, 8 and $s+1$, respectively;
\item all remaining  equicut polytopes - the snub Cube, the snub Dodecahedron
and $Aprism_{2s}$ $(s>1)$ - all have scale 2 and sizes 
$9/2$, $15/2$ and $s+1/2$, respectively;
\item the remaining, $Prism_{2s+1}$ $(s\geq 1)$ 
has scale 2 and size $s+3/2$; it is not an equicut graph. 
\end{enumerate}

Among all Catalan (dual Archimedean) $\l$-polyhedra:
\begin{enumerate}
\item all the zonohedra are as follows:
dual Cuboctahedron and dual Icosidodecahedron of sizes  4 and 6, respectively;
\item the only other doubling is dual truncated Icosahedron (also known as
dual football) of scale 2 and size  5;
\item all remaining cases (duals of truncated Cube, of truncated Dodecahedron, 
  of truncated Tetrahedron and of $Prism_3$) are non-equicut $\l$-graphs of
  scale 2 and of sizes 6, 13, $7/2$ and 2, respectively.
\end{enumerate}
All Archimedean or Catalan  $\l$-polyhedra are rigid 
except the dual $Prism_3$; those embeddings are presented in \cite{DeSt96}. 

Examples of regular-faced polyhedra (from the well-known list of 92
polyhedra, see \cite{Joh66,Za69}) with equicut skeletons, are (all of
scale 2) \#75 biaugmented $Prism_6$ of size 4 (a doubling) and two
non-doublings: \#74 augmented $Prism_6$ of size 4 and
\#83 tridiminished Icosahedron of size 3.

The regular $\l$-polytopes of dimension greater than 3 have 
equicut skeletons. They are as follows.
\begin{itemize}
\item $K_v$ ($(v-1)$-simplex);
\item $H(s,2)$ ($s$-cube);
\item $K_{s\times 2}$ (cross-polytope).
\end{itemize}

There are just 3 semiregular $\l$-polytopes of dimension greater than
3, see \cite{DeSh96}.  Two of them have equicut skeletons: 
$\frac{1}{2}H(5,2)$ and the snub 24-cell $s(3,4,3)$,
see Figure~\ref{fig:sn24}.  The latter is a
4-dimensional semiregular polytope with 96 vertices (see, for example,
\cite{Cox73}); the regular 4-polytope 600-cell can be obtained by
capping its 24 icosahedral facets.  Its skeleton has
scale 2 and size 6; it is a doubling.

Three of the chamfered (see \cite{DeSh96}) Platonic solids have $\l$-skeletons:
chamfered Cube is a zonohedron of size 7, chamfered Dodecahedron is an
equicut (non-doubling) graph of scale 2 and size 11, chamfered Tetrahedron
is non-equicut $\l$-graph of scale 2 and size 4.

\begin{figure}
\arraycolsep 0.38mm
\scriptsize
\[
\mathtt{
\begin {array}{cccccccccccccccccccccccccccccccccccccccccccccccc} 
0&0&0&0&0&0&0&0&0&0&0&0&0&0&0&0&0&0&0&0&0&0&0&0&0&0&0&0&0&0&0&0&0&0&0&0&0&0&0&0&0&0&0&0&0&0&0&0\\
0&0&0&0&0&0&0&0&0&0&0&0&0&0&0&0&0&0&0&0&0&0&0&0&1&1&1&1&1&1&1&1&1&1&1&1&1&1&1&1&1&1&1&1&1&1&1&1\\
0&0&0&0&0&0&0&0&0&0&0&0&1&1&1&1&1&1&1&1&1&1&1&1&0&0&0&0&0&0&0&0&0&0&0&0&0&0&0&0&0&0&0&0&1&1&1&1\\
0&0&0&0&1&1&1&1&1&1&1&1&0&0&0&0&0&0&0&0&1&1&1&1&0&0&0&0&0&0&0&0&1&1&1&1&1&1&1&1&1&1&1&1&0&0&0&0\\
0&0&0&0&0&0&0&0&1&1&1&1&0&0&0&0&1&1&1&1&1&1&1&1&0&0&0&0&1&1&1&1&0&0&0&0&0&0&0&0&1&1&1&1&1&1&1&1\\
0&0&0&0&0&0&0&0&0&0&0&0&0&0&0&1&0&0&1&1&0&0&0&1&0&0&1&1&0&1&1&1&0&0&0&0&0&0&0
&1&0&0&1&1&1&1&1&1\\
0&0&1&1&0&0&0&0&0&0&1&1&0&0&1&1&1&1&1&1&0&0&1&1&1&1&1&1&1&1&1&1&0&0&0&0&1&1&1&1&1&1&1&1&1&1&1&1\\
0&0&0&1&0&0&1&1&1&1&1&1&0&0&0&0&0&1&0&1&0&1&1&1&0&1&0&1&1&1&1&1&0&1&1&1&1&1&1&1&1&1&1&1&0&1&1&1\\
0&1&1&1&0&1&0&1&0&1&1&1&1&1&1&1&1&1&1&1&1&1&1&1&0&1&1&0&1&0&1&1&0&0&0&0&0&0&1&0&0&1&0&1&1&1&1&1\\
0&0&0&0&0&0&0&0&1&0&0&1&0&0&0&0&0&0&0&0&0&1&0&1&0&0&0&0&0&1&0&1&0&0&1&1&0&1&0&1&1&1&1&1&0&0&1&1\\
0&1&0&1&1&1&1&1&1&1&1&1&0&1&1&0&1&1&1&1&1&1&1&1&0&0
&0&0&1&0&0&1&0&1&0&1&0&1&1&0&1&1&1&1&0&1&0&1\\
0&0&0&0&0&1&1&0&1&1&0&1&0&1&0&0&1&0&0&1&1&1&1&1&0&0&0&0&0&0&0&0&0&0&0&1&0&0&0
&0&1&0&0&1&0&0&0&1
\end {array}}
\]
\caption{An equicut realization of the snub 24-cell; only half of the 96 column vectors
are shown. The remaining 48 are obtained as complements.}\label{fig:sn24}
\end{figure}

\medskip

{\em Authors' addresses}:
\medskip

\begin{tabular}{ll}
Michel Marie Deza          & Dmitrii V. Pasechnik\\
LIENS, DMI                 & Dept. of Computer Science\\
Ecole Normale Sup\'erieure~~~~~~~~~~ & Utrecht University\\
45 rue d'Ulm               &PO Box 80.089\\ 
75230 Paris Cedex 05       &3508 TB Utrecht\\
France                     &The Netherlands\\
{\em e-mail}: deza@dmi.ens.fr & {\em e-mail}: dima@cs.uu.nl
\end{tabular}
\end{document}